\documentclass[11pt]{amsart}
\usepackage{amsfonts}
\usepackage{amsmath}
\usepackage{amsthm}
\usepackage{amssymb}
\usepackage{latexsym}
\usepackage{multicol}
\usepackage{verbatim}

\advance\textwidth by 1.2in \advance\oddsidemargin by -.6in
\advance\evensidemargin by -.6in
\newtheorem*{cor}{Corollary}%[section]
\newtheorem*{lem}{Lemma}
\newtheorem*{prop}{Proposition}

\theoremstyle{definition}
\newtheorem*{defn}{Definition}
\theoremstyle{definition}
\newtheorem*{thm}{Theorem}

\newtheorem*{rem}{Remark}

\newenvironment{pf}{\proof}{\endproof}
\newcounter{cnt}
\newenvironment{enumerit}{\begin{list}{{\hfill\rm(\roman{cnt})\hfill}}{%
\settowidth{\labelwidth}{{\rm(iv)}}\leftmargin=\labelwidth%
\advance\leftmargin by
\labelsep\rightmargin=0pt\usecounter{cnt}}}{\end{list}}

\theoremstyle{remark}

\numberwithin{equation}{section} 

\def\wt{{\rm wt}}

\def\opl_#1^#2{\text{\scriptsize$\bigoplus\limits_{\text{\footnotesize$#1$}}^{\text{\footnotesize$#2$}}$}}

\begin{document}

\newcommand{\thmref}[1]{Theorem~\ref{#1}}
\newcommand{\secref}[1]{Section~\ref{#1}}
\newcommand{\lemref}[1]{Lemma~\ref{#1}}
\newcommand{\propref}[1]{Proposition~\ref{#1}}
\newcommand{\corref}[1]{Corollary~\ref{#1}}
\newcommand{\remref}[1]{Remark~\ref{#1}}
\newcommand{\defref}[1]{Definition~\ref{#1}}
\newcommand{\er}[1]{(\ref{#1})}
\newcommand{\id}{\operatorname{id}}
\newcommand{\tensor}{\otimes}
\newcommand{\nc}{\newcommand}
\newcommand{\rnc}{\renewcommand}
\newcommand{\qbinom}[2]{\genfrac[]{0pt}0{#1}{#2}}
\nc{\cal}{\mathcal} \nc{\goth}{\mathfrak} \rnc{\bold}{\mathbf}
\renewcommand{\frak}{\mathfrak}
\newcommand{\desc}{\operatorname{desc}}
\newcommand{\Maj}{\operatorname{Maj}}
\renewcommand{\Bbb}{\mathbb}
\nc\bomega{{\mbox{\boldmath $\omega$}}} \nc\bpsi{{\mbox{\boldmath
$\Psi$}}}
 \nc\balpha{{\mbox{\boldmath $\alpha$}}}
 \nc\bpi{{\mbox{\boldmath $\pi$}}}

\newcommand{\lie}[1]{\mathfrak{#1}}
\makeatletter
\def\section{\def\@secnumfont{\mdseries}\@startsection{section}{1}%
  \z@{.7\linespacing\@plus\linespacing}{.5\linespacing}%
  {\normalfont\scshape\centering}}
\def\subsection{\def\@secnumfont{\bfseries}\@startsection{subsection}{2}%
  {\parindent}{.5\linespacing\@plus.7\linespacing}{-.5em}%
  {\normalfont\bfseries}}
\makeatother
\def\subl#1{\subsection{}\label{#1}}
 \nc{\Hom}{\text{Hom}}
\nc{\ch}{\text{ch}}
\nc{\ev}{\text{ev}}

\nc{\krsm}{KR^\sigma(m\omega_i)}
\nc{\krsmzero}{KR^\sigma(m_0\omega_i)}
\nc{\krsmone}{KR^\sigma(m_1\omega_i)}
 \nc{\vsim}{v^\sigma_{i,m}}

 \nc{\Cal}{\cal} \nc{\Xp}[1]{X^+(#1)} \nc{\Xm}[1]{X^-(#1)}
\nc{\on}{\operatorname} \nc{\Z}{{\bold Z}}
\nc{\J}{{\cal J}} \nc{\C}{{\bold C}} \nc{\Q}{{\bold Q}}
\renewcommand{\P}{{\cal P}}
\nc{\N}{{\Bbb N}} \nc\boa{\bold a} \nc\bob{\bold b} \nc\boc{\bold c}
\nc\bod{\bold d} \nc\boe{\bold e} \nc\bof{\bold f} \nc\bog{\bold g}
\nc\boh{\bold h} \nc\boi{\bold i} \nc\boj{\bold j} \nc\bok{\bold k}
\nc\bol{\bold l} \nc\bom{\bold m} \nc\bon{\bold n} \nc\boo{\bold o}
\nc\bop{\bold p} \nc\boq{\bold q} \nc\bor{\bold r} \nc\bos{\bold s}
\nc\bou{\bold u} \nc\bov{\bold v} \nc\bow{\bold w} \nc\boz{\bold z}
\nc\boy{\bold y}
\nc\ba{\bold A} \nc\bb{\bold B} \nc\bc{\bold C} \nc\bd{\bold D}
\nc\be{\bold E} \nc\bg{\bold G} \nc\bh{\bold H} \nc\bi{\bold I}
\nc\bj{\bold J} \nc\bk{\bold K} \nc\bl{\bold L} \nc\bm{\bold M}
\nc\bn{\bold N} \nc\bo{\bold O} \nc\bp{\bold P} \nc\bq{\bold Q}
\nc\br{\bold R} \nc\bs{\bold S} \nc\bt{\bold T} \nc\bu{\bold U}
\nc\bv{\bold V} \nc\bw{\bold W} \nc\bz{\bold Z} \nc\bx{\bold x}
{\title[]{The restricted   Kirillov--Reshetikhin
modules for the current and twisted current algebras}

\author{Vyjayanthi Chari and  Adriano Moura}
\address{Department of Mathematics, University of
California, Riverside, CA 92521.}\thanks{VC was partially
supported by the NSF grant DMS-0500751}
\address{UNICAMP - IMECC, Campinas - SP - Brazil, 13083-859.} \email{aamoura@ime.unicamp.br}

\maketitle

\section*{Introduction}

In this paper we define and study a family of  $\bz_+$--graded
modules for  the polynomial valued current algebra $\lie g[t]$ and the twisted current algebra $\lie g[t]^\sigma$ associated
to a finite--dimensional classical simple Lie algebra $\lie g$ and a non--trivial diagram automorphism of $\lie g$.
The modules which we denote as $KR(m\omega_i)$ and $\krsm$ respectively
 are indexed by pairs $(i,m)$, where $i$ is a node  of the
Dynkin diagram and $m$ is a non--negative integer, and are given by
generators and relations. These modules are indecomposable, but  usually reducible, and we
describe their Jordan--Holder series by giving  the
corresponding graded  decomposition  as a direct sum of
irreducible modules for the underlying finite--dimensional simple
Lie algebra. Moreover, we
prove that the modules are finite--dimensional and hence restricted,
 i.e.,  there exists an integer $n\in\bz_+$
 depending only on $\lie g$ and $\sigma$  such that $(\lie g\otimes t^n)KR(m\omega_i)=0$.
It turns out that this graded decomposition is
exactly the one predicted in \cite[Appendix A]{HKOTY}, \cite[Section 6]{HKOTT}
coming from the study of the Bethe Ansatz in solvable lattice models.

Our interest in these modules
 and the motivation for calling them the Kirillov--Reshetikhin modules arises from the fact that when we
 specialize the grading by setting  $t=1$, (or equivalently putting $q=1$ in the formulae in \cite{HKOTY}, \cite{HKOTT})
the character of the module is exactly the one  predicted in \cite{KR},
 \cite{Kl} for a family of irreducible finite--dimensional modules for the Yangian of $\lie g$.
 Analogous modules for the untwisted  quantum affine algebra associated to $\lie g$ are also known
to exist with this decomposition \cite{Ckirres} and these have been studied from a combinatorial viewpoint in \cite{HKOTY},\cite{NS1}, \cite{NS2}, \cite{OSS}.
One of the methods used in \cite{Ckirres} involves passing to the $q=1$ limit of the modules for the quantum affine algebra,
although the resulting modules are not graded or restricted. The methods used
in \cite{Ckirres}  require a number of complicated results from the representation theory of the untwisted
 quantum affine algebras which have not been proved in the twisted case.

  In this paper we show that the Kirillov--Reshetikhin modules can be  studied in
the non--quantum case.
 Thus we define  graded restricted analogues  and compute their graded characters without resorting to the quantum situation. As  a consequence we have a mathematical interpretation of the parameter $q$ which appears in the fermionic formulae in \cite{HKOTY}, \cite{HKOTT}.  Then we prove that the abstractly defined modules $KR(m\omega_i)$ have a concrete construction as follows:  the $\lie g[t]$--module structure of the  \lq\lq fundamental\rq\rq \ Kirillov--Reshetihkhin modules (in most cases these modules correspond to taking $m=1$) is described explicitly and then  the module $KR(m\omega_i)$ is realized  as a canonical submodule of a tensor product of the fundamental modules.

 Another motivation for our interest is the connection with Demazure modules, \cite{CL},\cite{FL},\cite{FFK}, \cite{FoL}.
 Thus we are able to prove   that the Kirillov--Reshetikhin modules for the current algebras  are isomorphic
as representations of the current algebra  to the Demazure modules
in multiples of the basic representation of the current algebras.
Our methods also work for some nodes of the other exceptional
algebras, we explain this together with the reasons for the
difficulties for the  exceptional algberas in the concluding
section of this paper.

\setcounter{section}{0}
\section{Preliminaries}

\subsection{Notation}  Let $\bz_+$ (resp.  $\bn$) be the set of
non-negative (resp.  positive) integers. Given  a Lie algebra
$\lie a$, let  $\bu(\lie a)$ denote  the universal enveloping  algebra of
$\lie a$ and  $\lie a[t]=\lie a\otimes \bc[t]$  the
polynomial valued current  Lie algebra of $\lie a$. The Lie algebra
$\lie a[t]$ and its universal enveloping algebra  are  $\bz_+$-graded where the grading is given by the
powers of $t$. We shall identify $\lie a$ with the subalgebra $\lie a\otimes 1$  of $\lie a[t]$.

 \subsection{Classical simple Lie algebras} For the rest of the paper $\frak g$
 denotes  a complex finite--dimensional simple
Lie algebra   of type
$A_n$, $B_n$, $C_n$ or $D_n$, $n\ge 1$, and $\frak h$  a  Cartan
subalgebra of $\lie g$. Let $I=\{1,\cdots, n\}$ and
    $\{\alpha_i:i\in I\}$ (resp.  $\{\omega_i:i\in I\}$)
  be a set of simple roots
 (resp.  fundamental weights) of $\frak g$ with respect to $\frak h$,
   $R^+$ (resp. $Q$, $P$) be the corresponding set of positive roots (resp. root lattice, weight lattice) and let $Q^+$, $P^+$ be the $\bz_+$--span of the simple roots and fundamental weights respectively. It is convenient to set $\omega_0=0$.  Let $\theta\in R^+$
be the highest root. For  $i\in I$,  let $\varepsilon_i: Q\to \bz$ be defined by requiring $\eta=\sum_{i\in I}\varepsilon_i(\eta)\alpha_i$, $\eta\in Q$.
We assume throughout the paper that the simple roots are indexed as in \cite{Bo}.
Given $\alpha\in R$, let $\lie g_\alpha$ be the corresponding root
space. Fix non--zero  elements $x^\pm_\alpha\in\lie g_{\pm\alpha}$,
$h_\alpha\in\lie h$, such that $$[h_\alpha,x^\pm_\alpha]=\pm
2x^\pm_\alpha,\ \ \ \ [x^+_\alpha,x^-_\alpha]=h_\alpha.$$
 Set $\lie n^\pm=\opl_{\alpha\in R^+}^{}\lie g_{\pm\alpha}$.
Given any subset $J\subset I$ let $\lie g(J)$ be the subalgebra of
$\lie g$ generated by the elements $x^\pm_{\alpha_j}$, $j\in
J$. The sets $R(J)$, $P(J)$ etc. are defined in the obvious way.

Let $(\ ,\ )$ be the form on $\lie h^*$ induced by the restriction
of the  Killing form of $\lie g$ to $\lie h$ normalized so that $(\theta,\theta)=2$  and set
$\check{d_j}=2/(\alpha_j,\alpha_j)$.

\subsection{Finite--dimensional $\lie g$--modules}

Given $\lambda\in P^+$, let $V(\lambda)$ be the irreducible
finite--dimensional $\lie g$--module with highest weight vector
$v_\lambda$, i.e., the cyclic module generated by $v_\lambda$ with
defining relations: $$\lie n^+v_\lambda=0,\ \
hv_\lambda=\lambda(h)v_\lambda,\ \
(x^-_{\alpha_i})^{\lambda(h_{\alpha_i})+1}v_\lambda=0.$$ For $\mu\in
P$, let $V(\lambda)_\mu=\{v\in V(\lambda): hv=\mu(h)v,\ \ h\in \lie
h\}. $   Any finite--dimensional $\lie
g$--module $M$ is  isomorphic to a direct sum $\oplus_{\lambda\in
P^+}V(\lambda)^{\oplus m_\lambda(M)}$, $m_{\lambda}(M)\in\bz_+$.

\begin{prop}\cite{PRV}\label{tpgen}  Let $\lambda,\mu\in P^+$. Then $V(\lambda)\otimes
V(\mu)$ is generated as a $\lie g$--module by the element
$v_\lambda\otimes v_\mu^*$ and relations:
$$h(v_\lambda\otimes v^*_\mu)=(\lambda+\mu^*)(h)(v_\lambda\otimes
v_\mu^*),$$ and
$$(x^+_{\alpha})^{-\mu^*(h_\alpha)+1}(v_\lambda\otimes v^*_\mu)=(x^-_{\alpha})^{\lambda(h_\alpha)+1}(v_\lambda\otimes v^*_\mu)=0,$$
for all $\alpha\in R^+$. Here $\mu^*$  is  the lowest weight of $V(\mu)$ and $0\ne v_\mu^*\in V(\mu)_{\mu^*}$.\hfill\qedsymbol
\end{prop}

\subsection{Graded Modules}

Given a $\lie g$--module $M$, we regard it as a  $\lie g[t]$--module by setting $(x\otimes t^r) m=0$ for all $m\in M$, $r\in\bn$, and denote the resulting
 graded  $\lie g[t]$--module by $\ev_0(M)$.
Let $V=\oplus_{s\in\bz_+}V[s]$ be a graded representation of $\lie g[t]$ with $\dim(V[s])<\infty$. Note that
each $V[s]$ is a $\lie g$--module.
For any $s\in\bz_+$, let $V(s)$ be the $\lie g[t]$--quotient of $V$ by the submodule $\oplus_{s'\ge s}V[s']$.
Clearly $$V(s)/V(s+1)\cong \ev_0(V[s]),$$  and  the irreducible constituents of $V$ are just the irreducible constituents of $
\ev_0(V[s])$, $s\in\bz_+$.

\section{The Kirillov-Reshetikhin modules for $\lie g[t]$}

 \subsection{The modules $KR(m\omega_i)$} \begin{defn} Given $i\in I$  and $m\in\bz_+$,
let $KR(m\omega_i)$ be the $\lie g[t]$--module generated by an
element $v_{i,m}$ with relations, \begin{equation}\label{defrel}\lie n^+[t] v_{i,m}=0,\ \
h v_{i,m}=m\omega_i(h)v_{i,m},\ \ (h\otimes t^r)v_{i,m}=0,\ \ \
h\in\lie h, r\in \bn, \end{equation} and \begin{equation}\label{defrel2}(x^-_{\alpha_i})^{m+1}v_{i,m}=
(x^-_{\alpha_i}\otimes t)v_{i,m}=x^-_{\alpha_j}v_{i,m}=0 ,\ \ j\in
I\backslash\{ i\}.\end{equation} \hfill\qedsymbol\end{defn} Note that the modules $KR(m\omega_i)$ are graded
modules since the defining relations are graded.
The following is trivially checked.
\begin{lem}For all $i\in I$, $m\in\bz_+$, the module $\ev_0(V(m\omega_i))$
is a quotient of $KR(m\omega_i)$. \hfill\qedsymbol\end{lem}

\begin{rem} These modules were defined and  studied initially in  \cite[Section 1]{Ckirres} (where
they were denoted as $W(i,m)$) as quotients of the finite--dimensional
Weyl modules defined in \cite{CPweyl}.
In this paper however, we shall show
directly  that the modules  $KR(m\omega_i)$ are finite--dimensional
  as a consequence of the analysis of
their $\lie g$--module structure.
\end{rem}

\subsection{The graded character of $KR(m\omega_i)$} We now state the main result of this section. We need some additional notation. For $i\in I$, $m\in\bz_+$, let $P^+(i,m)
\subset P^+$ be defined  by:
$$P^+(i,m)= P^+(i,\check d_i)+P^+(i,m-\check d_i),\ \ m> \check d_i,$$ where $P^+(i,1)=\omega_i$ if
  $\varepsilon_i(\theta)=\check d_i$ and in the other cases,
\begin{eqnarray*} P^+(i,1)&=&\{\omega_i,\omega_{i-2},\cdots
,\omega_{\overline i}\},\  \ \ \ \  \lie g= D_n, \\
 P^+(i,1)&=&\{\omega_i,\omega_{i-2},\cdots
,\omega_{\overline i}\},\ \ \ \  \lie g=B_n,\\
P^+(n,2)&=&\{2\omega_n,\omega_{n-2},\cdots ,\omega_{\overline n}\},
\ \ \lie g=B_n,\\ P^+(i,2)&=&\{2\omega_i,2\omega_{i-1},\cdots
,2\omega_1, 0\},\ \ \  \ \ \lie g=C_n,
 \end{eqnarray*}
where $ \overline i\in\{0,1\}$ and $i=\overline i\mod 2$. Let  $\mu_0,\cdots ,\mu_k$ be the unique enumeration of the sets $P^+(i,\check d_i)$
chosen so  that $$\mu_j-\mu_{j+1}\in R^+,\ \ \mu_j-\mu_{j+2}\in Q^+\backslash R^+,\ \ \ 0\le j\le k-1.$$ Given $m=\check d_im_0+m_1$ with $0\le m_1<\check d_i$ and
$\mu\in P^+(i,m)$, we can clearly write $\mu=m_1\omega_i+\mu_{j_1}+\cdots +\mu_{j_{m_0}}$,
where $j_r\in\{0,1,\cdots ,k\}$ for $1\le r\le m_0$. We say that the expression is reduced if each
 $j_r$ is minimal with the property that $\mu-\mu_{j_1}-\cdots -\mu_{j_r}\in P^+(i,m-r)$.  Such an expression is clearly unique and we set $$|\mu|=\sum_{r=1}^{m_0}j_r.$$

\begin{thm}\label{mainun}\hfill
\begin{enumerit}
\item Let $i\in I$, $m,s\in\bz_+$. We have,
$$KR(m\omega_i)[s]\cong\bigoplus_{\{\mu\in P^+(i,m):|\mu|=s\}} V(\mu).$$

\item  Write $m=\check d_im_0+m_1$, where $0\le
m_1<\check d_i$. The canonical homomorphism of $\lie g[t]$--modules  $$
KR(m\omega_i)\to KR(m_1\omega_i)\otimes KR(\check d_i\omega_i)^{\otimes
m_0}$$ mapping $v_{i,m}\mapsto v_{i,m_1}\otimes v_{i,\check d_i}^{\otimes m_0}$ is injective.
\end{enumerit}
\end{thm}

The rest of the section is devoted to the proof of this result.

 \subsection{Elementary properties of $KR(m\omega_i)$} \begin{prop}\label{eleprop}\hfill
 \begin{enumerit}
\iffalse\item $KR(0)\cong \bc$. \item For all $\alpha\in R^+$ and
$s\in\bz_+$, we have
$$(x^-_\alpha\otimes t^s)^{m\omega_i(h_\alpha)+1}v_{i,m}=0.$$\fi \item
We have
$$KR(m\omega_i)=\opl_{\mu\in\lie h^*}^{} KR(m\omega_i)_\mu $$ and $KR(m\omega_i)_\mu\ne 0$ only if
$\mu\in m\omega_i-Q^+$.
\item Regarded as a $\lie g$--module, $KR(m\omega_i)$ and $KR(m\omega_i)[s]$, $s\in\bz_+$,
are isomorphic to a direct sum of irreducible finite--dimensional
representations of $\lie g$.
\item For all $0\le r\le m$,  there exists a canonical homomorphism
$KR(m\omega_i)\to KR(r\omega_i)\otimes
KR((m-r)\omega_i)$ of graded $\lie g[t]$--modules such that
$v_{i,m}\mapsto v_{i,r}\otimes v_{i,(m-r)}$.
\end{enumerit}
\end{prop}

\begin{pf} Part (i) follows by a
standard application of the PBW theorem. For (ii) it suffices, by
standard results, to show that $KR(m\omega_i)$ is a sum of
finite--dimensional $\lie g$--modules. Note  that the defining
relations of $KR(m\omega_i)$ imply that  $\bu(\lie g)v_{i,m}\cong
V(m\omega_i)$ as  $\lie g$--modules. Hence
\begin{equation} \label{locnil} (x^-_\alpha)^{m\omega_i(h_\alpha)+1}v_{i,m}=0\end{equation} for all
$\alpha\in R^+$. \iffalse For all $s\in\bn$ we have
$$(h_\alpha\otimes
t^{s})(x^-_\alpha)^{m\omega_i(h_\alpha)+1}v_{i,m}=
(x^-_\alpha)^{m\omega_i(h_\alpha)}(x^-_\alpha\otimes t^s)v_{i,m}=
0.$$ Repeating this we get (i) and (ii). \fi For $v\in
KR(m\omega_i)_\mu$ we have  $\bu(\lie g)v=\bu(\lie n^-)\bu(\lie
n^+)v$. Part  (ii) implies that $\bu(\lie n^+)v$ is a
finite--dimensional vector space. Further since the  action of
$\lie n^-$ on $\lie n^-[t]$ given by the Lie bracket is locally
nilpotent and since $v\in\bu(\lie n^-[t])v_{i,m}$, it follows from
\eqref{locnil}  that $x^-_\alpha$ acts nilpotently on $v$.  This
proves that $\bu(\lie g)v$ is finite--dimensional and part (ii) is
established. Part (iii) is clear from the defining relations of
the modules.
\end{pf}

\begin{cor} For $i\in I$, $m\in\bz_+$,
we have $$KR(m\omega_i)=\opl_{\mu\in P^+}^{}V(\mu)^{\oplus
m_\mu(i,m)}.$$ In particular,
$KR(0)\cong\bc$.\hfill\qedsymbol\end{cor}

\subsection{An upper bound for $m_\mu(i,m)$.} The following result was proved in \cite[Theorem 1
]{Ckirres} under the assumption that $KR(m\omega_i)$ is
finite--dimensional. An inspection of the proof shows however,
that the only place this is used is to write $KR(m\omega_i)$ as a
direct sum of irreducible $\lie g$--modules. But (as we shall  see in the twisted case, where we do give a proof of the analogous proposition)  this only
requires the weaker result proved in Proposition \ref{eleprop},
 \begin{prop}\label{unupb}  As a $\lie g$-module we have $$KR(m\omega_i)\cong
 \opl_{\mu\in P^+(i,m)}^{}V(\mu_{})^{\oplus
m_{\mu}},$$ where $m_{\mu}\in\{0,1\}$.
\end{prop}
Th next corollary is  immediate since $P^+(i,m)$ is a finite  set by definition.
\begin{cor} For all $i\in I$,
the modules $KR(m\omega_i)$ are finite--dimensional.\hfill\qedsymbol\end{cor}

\subsection{ Proof of Theorem \ref{mainun}: the cases $\varepsilon_i(\theta)=1$, $m\in\bz_+$ and
 $\varepsilon_i(\theta)=\check d_i$, $m=1$}\hfill
In these cases, it follows from the definition that $P^+(i,m)=\{m\omega_i\}$.
Since $\ev_0(V(m\omega_i))$ is a $\lie g[t]$--module quotient of $KR(m\omega_i)$, part (i) is immediate
 from Proposition \ref{unupb}.
Part (ii) of the theorem is now  obvious  since the canonical inclusion
$V(m\omega_i)\to V(\omega_i)^{\otimes m}$ of $\lie
g$--modules is obviously also an inclusion of the  $\lie
g[t]$--modules $\ev_0(V(m\omega_i))\to\ev_0(V(\omega_i))^{\otimes m}$.

\vskip 24pt

{\it{ Since $\varepsilon_i(\theta)=1$ for all $1\le i\le n$ if $\lie g$ is of type $A_n$, the theorem is proved in this case and  we assume for the rest of this section   that $\lie g$ is not of type $A_n$.}}

\subsection{An explicit construction in the case $\varepsilon_i(\theta)=2$ and $m=\check d_i$.} \hfill

Let $V_s$,  $0\le s\le k$,  be  $\lie g$--modules  such that
\begin{equation}\label{condun}  \Hom_{\lie g}(\lie g\otimes V_s,
 V_{s+1})\ne 0,\ \ \Hom_{\lie g}(\wedge^2(\lie g)\otimes V_s,
 V_{s+2})=0,\ \ 0\le s\le k-1,\end{equation} where we assume that $V_{k+1}=0$. Fix non--zero elements $p_s\in\Hom_{\lie g}(\lie g\otimes V_s, V_{s+1})$, $0\le s\le k-1$, and set $p_k=0$.

It is easily checked that the following formulas extend the canonical  $\lie g$--module structure to   a graded  $\lie
 g[t]$--module structure on $V=\oplus_{s=1}^k V_s$:
 $$ (x\otimes t)v=p_s(x\otimes v), \ \
 (x\otimes t^r)v=0,\ \ r\ge 2,$$ for all $
 x\in\lie g$, $v\in V_s,  1\le s\le k$.
Clearly, $$V[s]
\cong _{\lie
  g}V_{s},\ \ 0\le s\le k.$$
Moreover, if the maps $p_s$, $0\le s\le k-1$, are all surjective and if  $V_0=\bu(\lie g)v_0$
 then $V=\bu(\lie g[t])v_0$.

\begin{prop}\label{constun} Let $i\in I$ be such that $\varepsilon_i(\theta)=2$ and let $\mu_s\in P^+(i,\check d_i)$, $0\le s\le k$.
  The modules
 $V(\mu_s)$, $0\le s\le k$, satisfy \eqref{condun} and the resulting $\lie g[t]$--module is isomorphic to $KR(\check d_i\omega_i)$. In particular,
$$KR(\check d_i\omega_i)[j]\cong _{\lie g}V(\mu_j),\ \ 0\le j\le k.$$
\end{prop}
 \begin{pf}

 Using Proposition \ref{tpgen} it is easy to see that $$\Hom_{\lie g}(\lie g\otimes V(\mu_s), V(\mu_{s+1}))\cong \Hom_{\lie
 g}(\lie g, V(\mu_s)\otimes V(\mu_{s+1}))\ne 0$$ for all $0\le s
 \le k-1$.
It is not hard to check (see \cite{FH}, \cite{V}) that as $\lie g$--modules,
 $$\wedge^2(\lie
 g)\cong\lie g\oplus V(\nu),$$ where $\nu=2\omega_1+\omega_2$ if
 $\lie g$ is of type $C_n$, $\nu=\omega_1+2\omega_3$ if $\lie g$ is of type $B_3$, $\nu=\omega_1+\omega_3+\omega_4$ if $\lie g$ is of type $D_4$, and $\nu=\omega_1+\omega_3$ otherwise.
 Since $\mu_s-\mu_{s+2}\notin R$, it follows  that
 $\Hom_{\lie g}(\lie g\otimes V(\mu_s), V(\mu_{s+2}))=0$.
To prove that $\Hom_{\lie g}(V(\nu)\otimes V(\mu_s),
V(\mu_{s+2}))=0$, it suffices to prove that   $$\Hom_{\lie
g}(V(\mu_s), V(\nu)\otimes V(\mu_{s+2}))=0.$$ Suppose that $0\ne
p\in \Hom_{\lie g}( V(\mu_{s}),V(\nu)\otimes V(\mu_{s+2}))$. A simple computation using the explicit formulas for the fundamental weights in terms of the simple roots (see \cite{H} for instance) shows that
$$\mu_{s+2}+\nu-\mu_s\in Q^+_J,\ \  J=\{1,\cdots, s-1\}.$$ Hence
 $$0\ne p(v_{\mu_s})\in \bu(\lie
g_J)v_\nu\otimes\bu(\lie g_J)v_{\mu_{s+2}},$$ which implies that
$$p(\bu(\lie g_J)v_{\mu_s})\subset
\bu(\lie g_J)v_\nu\otimes\bu(\lie g_J)v_{\mu_{s+2}}.$$  The formulas given for $\mu_s$ shows that  $\mu_s(\lie h_J)=0$. This means that $\bu(\lie g_J)v_{\mu_s}\cong\bc$ which implies that
$$p(\bc)\cong\bc\subset \left( \bu(\lie g_J)v_\nu\otimes\bu(\lie g_J)v_{\mu_{s+2}}\right).$$
But this is impossible since by standard results $\bu(\lie
g_J)v_{\mu_{s+2}}$ is not isomorphic to the $\lie g_J$--dual of $\bu(\lie g_J)v_{\nu} $. Hence $p=0$ and
 we have proved that the modules $V(\mu_s)$, $0\le s\le k-1$, satisfy the conditions of \eqref{condun}. Let $\tilde{KR}(\check d_i\omega_i)$ be the resulting $\lie g[t]$--module. Since the modules $V(\mu_s)$ are all irreducible, it follows moreover that $$\tilde{KR}(\check d_i\omega_i)=\bu(\lie g[t])v_{\check d_i\omega_i}.$$

To see that $\tilde{KR}(\check d_i\omega_i)\cong KR(\check d_i\omega_i)
$ it suffices, by Proposition \ref{unupb}, to prove that $\tilde{KR}(\check d_i\omega_i)$ is a quotient of $KR(\check d_i\omega_i)$.  Since  $\mu_0=m\omega_i$ and  $\mu_0-\mu_1\in
R^+\backslash\{\alpha_i\}$, we see that we must have,$$p_0\left((\lie n^+\oplus
\lie h\oplus\bc x^-_{\alpha_i})\otimes v_{m\omega_i}\right)=0.$$
This proves that $$\lie n^+[t]v_{m\omega_i}=0, \ \ h
v_{m\omega_i}=(m\omega_i)(h) v_{m\omega_i}, \ \ (h\otimes
t^r)v_{m\omega_i}= \ \ (x^-_{\alpha_i}\otimes t)v_{m\omega_i}=0,\
\ \ r\ge 1.$$ Finally, since $\tilde{KR}(m\omega_i)$ is obviously finite--dimensional, it follows that $$(x^-_{\alpha_i})^{m\omega_i(h_{\alpha_i})+1}v_{\mu_0}=x^-_{\alpha_j}v_{\mu_0}=0,\ \ \ j\in I\backslash\{i\},$$ and the proof of the proposition is complete.
\end{pf}
Given $\mu\in P^+(i,m)$, $m=\check d_im_0+m_1$, $0\le m_1< \check d_i$,  with reduced expression $$\mu=m_1\omega_i+\mu_{j_1}+\cdots +\mu_{j_{m_0}},$$ set $s_j=\#\{r:j_r\ge j\}$, $0\le j\le k$, and  $$\bx_\mu=(x^-_{\mu_{k-1}-\mu_k}\otimes t)^{s_k}\cdots (x^-_{\mu_0-\mu_1}\otimes t)^{s_1}\in\bu(\lie n^-[t]).$$
It is easily  seen that  $$\bx_\mu v_{i,m}\in KR(m\omega_i)[|\mu|]\cap KR(m\omega_i)_{\mu}.$$

\begin{cor} \label{fundun}
Let $i\in I$ be such that $\varepsilon_i(\theta)=2$.
For  $0\le s\le k$ we have
$$\bx_{\mu_s} v_{i,\check d_i}\ne 0,\ \ h(\bx_{\mu_s}
v_{i,\check d_i})=\mu_s(h) (\bx_{\mu_s}v_{i,\check d_i}),\ \ \lie
n^+(\bx_{\mu} v_{i,\check d_i})=0, \ \ h\in\lie h,$$ or, equivalently, $$\bx_{\mu_s}v_{\mu_0}=v_{\mu_s}.$$

\end{cor}
\begin{pf} Suppose that $\mu=\mu_s$ for some $1\le s\le k$.
  We first prove that  there exists $x_s\in\lie n^-$ such that  \begin{equation}\label{req}p_s(x_s\otimes v_{\mu_s})=v_{\mu_{s+1}}.\end{equation} For this, note that by Proposition \ref{tpgen} $$v=p_s(x^-_\theta\otimes v_{\mu_s})\ne 0.$$
 Since $V(\mu_{s+1})$ is irreducible there exists
 $x_s'\in\bu(\lie n^+)$ such that   $x_s'v=v_{\mu_{s+1}}$, which gives $$x_s'(p_s(x^-_\theta\otimes v_{\mu_s}))=p_s({\text{ad}}(x_s')x^-_\theta\otimes v_{\mu_s})=v_{\mu_{s+1}},$$ proving \eqref{req}. Setting  $ x_s={\text{ad}}(x_s')x^-_\theta\in
\lie n^-_{\mu_s-\mu_{s+1}}$ it follows that $x_s$ is a non--zero scalar multiple of $x^-_{\mu_{s}-\mu_{s+1}}$ for some $0\le s\le k$
and we get $$(x_{\mu_{s}-\mu_{s+1}}^-\otimes t)v_{\mu_s}\ne 0.$$ An obvious induction on $s$ now proves that $$\bx_{\mu_s}v_{i,\check d_i}\ne 0,\ \ 1\le s\le k-1.$$ The other two statements of the corollary are now immediate.
\end{pf}
\vskip 12pt

Recall that $$KR(
\check d_i\omega_i)(r)\cong KR(\check d_i\omega_i)/\left(\oplus_{r'\ge r}KR(m\omega_i)[r']\right).$$  Let  $v_r$ be the image of $v_{\check{d_i}\omega_i}$ in $KR(\check d_i\omega_i)(r)$. By Corollary
   \ref{constun} we see that
   \begin{equation}\label{graded}\bx_{\mu_{r-1}}v_r\ne 0,\ \ \bx_{\mu_{r'}}v_r=0,\ \ 0\le r\le r'\le k.\end{equation}

\subsection{Proof of Theorem \ref{mainun}: the case  $\varepsilon_i(\theta)=2$, $m\ge \check d_i$} \hfill

If $m=\check d_i$, then part (i) is the statement of Proposition \ref{constun} and part (ii) is trivially true. Assume that $m>\check d_i$ and write $m=\check{d_i}m_0+m_1$,
$0\le m_1<\check{d_i}$. Note that since $m_1\in\{0,1\}$ we know by the earlier case that $$KR(m_1\omega_i)\cong \ev_0(V(m_1\omega_i)). $$
 Set
$$\tilde{KR}(m\omega_i)=\bu(\lie g[t])(v_{m_1\omega_i}\otimes v_{\check d_i\omega_i}^{\otimes m_0})\subset
KR(m_1\omega_i) \otimes KR(\check d_i\omega_i)^{\otimes m_0}.$$
 It is easily checked that
 $\tilde
{KR}(m\omega_i)$ is a graded quotient of $KR(m\omega_i)$. Using Proposition \ref{unupb} we see that part (i) of the theorem follows if we prove that for all $s\in\bz_+$,
\begin{equation}
\label{chtilde}\tilde{KR}(m\omega_i)[s]\cong \bigoplus_{\{\mu\in P^+(i,m):|\mu|=s\}}
V(\mu).\end{equation} In particular this proves that
$\tilde{KR}(m\omega_i)\cong KR(m\omega_i)$ and so
  proves part (ii) of the theorem as well.

Let $\mu=m_1\omega_i+\mu_{j_1}+\cdots+\mu_{j_{m_0}}$ be a reduced expression for $\mu$, set $$K= KR(m_1\omega_i)\otimes KR(\check d_i\omega_i)({j_1}+1)
\otimes\cdots\otimes KR(\check d_i\omega_i)({j_{m_0}}+1)
 ,$$
and let $\tilde K=\bu(\lie g[t])(v_{m_1\omega_i}\otimes v_{j_1}\otimes \cdots \otimes v_{j_{m_0}}).$
  Clearly $K$ and $\tilde K$ are  graded  quotients  of $KR(m_1\omega_i)\otimes  KR(\check d_i\omega_i)^{\otimes m}$ and $\tilde{KR}(m\omega_i)$ respectively.
 Let $$\pi:KR(m_1\omega_i)\otimes  KR(\check d_i\omega_i)^{\otimes
 m_0}\to K,\ \  \pi(\tilde{KR}(m\omega_i))= \tilde K,$$ be the canonical surjective morphism of $\lie g[t]$--modules.

Using \eqref{graded} and the comultiplication of $\bu(\lie g[t])$ one computes easily that
\begin{equation}\label{tilde2}\bx_\mu\left(\pi(v_{m_1\omega_i}
\otimes v_{\check d_i\omega_i}^{\otimes m_0})\right)=v_{m_1\omega_i}\otimes \bx_{\mu_{j_1}}v_{j_1}\otimes \cdots\otimes \bx_{\mu_{j_{m_0}}}v_{j_{m_0}}\ne 0.
\end{equation}

 Corollary \ref{fundun} implies that
$$\ \ \lie n^+\left(v_{m_1\omega_i}\otimes \bx_{\mu_{i_1}}v_{j_1}\otimes \cdots\otimes \bx_{\mu_{j_{m_0}}}v_{j_{m_0}}\right)=0,$$
and $$v_{m_1\omega_i}\otimes \bx_{\mu_{j_1}}v_{j_1}\otimes \cdots\otimes \bx_{\mu_{j_{m_0}}}
v_{j_{m_0}}\in
 K[|\mu|]\cap K_\mu.$$
It follows immediately that $$\tilde K[|\mu|]\cong_{\lie g} (V(\mu)\oplus N)$$ for some $\lie g$--submodule $N\subset\tilde K$ which proves \eqref{chtilde}.

\section{The twisted algebras }
We  use the notation of the previous section freely.

\subsection{Preliminaries} Throughout this section we let  $\epsilon \in\{0,1\}$.
From now on  $\lie g$ denotes a Lie algebra of type $A_n$ or $D_n$
and $\sigma:\lie g\to \lie g$ the non--trivial diagram
automorphism of order two. The statements in this subsection can be found in \cite{K}. Thus we have
$$\lie g=\lie g_0\oplus\lie g_1,\ \ \lie h=\lie h_0\oplus \lie h_1,\ \ \lie n^\pm =\lie
n_0^\pm \oplus \lie n_1^\pm ,$$ where $$\lie g_0=\{x\in\lie g:
\sigma(x)= x\},\ \  \lie g_1=\{x\in\lie g: \sigma(x)= -x\}.$$
For any subalgebra $\lie a$ of $\lie g$ with $\sigma(\lie a)\subset \lie a$ we set $\lie a_\epsilon=\lie a\cap\lie g_\epsilon$ and we have $\lie a=\lie a_0\oplus\lie a_1$.
The subalgebra $\lie g_0$
is a simple Lie algebra with Cartan subalgebra $\lie h_0$ and we
let $I_0$ be the  index set for the corresponding set of simpe roots numbered as in \cite{Bo}.

{\it{Although in this and the following sections it is convenient to use the ambient Lie algebra $\lie g$, we do not need any other data associated with it. Thus, all representations, roots, weights, the maps $\varepsilon_i$ and so on  will always be those associated with the fixed point algebra $\lie g_0$.}}

 We have $\lie g_0$
is of type $C_n$ if $\lie g$ is of type $A_{2n-1}$ and of type
$B_n$ if $\lie g $ is of type $A_{2n}$ or $D_{n+1}$. Let $(R_0)_s$
be the set of  short roots of $\lie g_0$ and $(R_0)_\ell$ the set
of long roots.
The adjoint action of $\lie g_0$ on $\lie g$ makes  $\lie g_1$ into an
irreducible representation of $\lie g_0$  and we have
 $$\lie g_1=\lie h_1\opl_{\mu\in\lie h_0^*}^{}(\lie g_1)_\mu,\ \ \lie g_1\cong V(\phi),$$
 where $\phi\in Q_0^+$ is the highest short root of $\lie g_0$ if $\lie g$
is of type $A_{2n-1}$ or $D_{n+1}$ and twice the highest short
root if $\lie g$ is of type $A_{2n}$.
 Further, if we set
$$R_1=\{\mu\in\lie h_0^*: (\lie g_1)_\mu\ne 0 \}\backslash\{0\},\ \
     R_1^+ = R_1\cap Q_0^+$$ then
$R_1=(R_0)_s$ if  $\lie g$  is of type $A_{2n-1}$ or $D_{n+1}$ and
$$R_1= R_0\cup\{ 2\alpha:\alpha\in (R_0)_s\}$$ if $\lie g$ is
of type $A_{2n}.$ In all cases, $\dim(\lie g_1)_\alpha=1$ for
$\alpha\in R_1$.

{\it{ Given $\alpha\in R_1^+$, we denote by
$y^\pm_{\alpha}$ {\bf {  any}} non--zero element in $(\lie g_1)_{\pm\alpha}$.
Note that $\lie n^\pm_1=\oplus_{\alpha\in R_1^\pm}^{} (\lie
g_1)_\alpha.$ In addition, if $\alpha\in R^+_0\backslash R^+_1$ we set $y_\alpha^\pm=0$, and similarly we set $x^\pm_\alpha=0$ if $\alpha\in R^+_1\backslash R_0^+$. }}

\begin{lem}\label{tbasic}\hfill
\begin{enumerit}
\item  The maps $\lie g_1\otimes \lie g_1\to \lie g_0$ and $\lie g_0\otimes \lie g_1\to \lie g_1$
defined by $x\otimes y\to [x,y]$ are surjective homomorphism of
$\lie g_0$--modules.

\item
 If $\alpha\in R_0\cap R_1$, there
exists $h\in\lie h_1$ such that $[h,y_\alpha^-] =
x^-_\alpha$.\end{enumerit}\hfill\qedsymbol
\end{lem}

\subsection{The twisted current algebra} Extend $\sigma$ to an automorphism $\sigma_t:\lie
g[t]\to\lie g[t]$ by extending linearly the assignment,
$$\sigma_t(x\otimes t^r)=\sigma(x)\otimes (-t)^r,\ \ x\in\lie g, \
\ r\in\bz_+.$$
 If $\lie a \subset\lie g$ is such that $\sigma(\lie a)\subset \lie a$, let $\lie a[t]^\sigma$ be the set of fixed points of $\sigma_t$. Clearly $$\lie a[t]^\sigma=\lie a_0\otimes \bc[t^2]\oplus \lie a_1\otimes t\bc[t^2].$$
and  $$\lie g[t]^\sigma=\lie
n^-[t]^\sigma\oplus\lie h[t]^\sigma\oplus \lie n^+[t]^\sigma. $$

\subsection{The modules $\krsm$} \hfill

\begin{defn} For  $i\in I_0$, $m\in\bz_+$, let $KR^\sigma(m\omega_i)$ be
the $\lie g[t]^\sigma$--module generated by an element
$v^\sigma_{i,m}$ with relations, \begin{equation}\label{deftw1}\lie n^+[t]^\sigma
v^\sigma_{i,m}=0,\ \ h_0v^\sigma_{i,m}=m\omega_i(h_0)v^\sigma_{i,m},\
\ (h_\epsilon \otimes t^{2r-\epsilon}) v_{i,m}^\sigma =x^-_{\alpha_j}v^\sigma_{i,m}=0,
\end{equation}
for all $h_\epsilon\in\lie h_\epsilon, r\in \bz_+, j\neq i,$
\begin{equation}\label{deftw2}(x^-_{\alpha_i})^{m+1}v^\sigma_{i,m}=0,
\end{equation} and
 \begin{equation}\label{deftw4} (x^-_{\alpha_i}\otimes t^2)\vsim
=(y_{\alpha_i}^-\otimes t)v^\sigma_{i,m}=0.\end{equation}
\hfill\qedsymbol\end{defn}
Note that when $\lie g$ is of type $A_{2n}$ or if $\alpha_i\in (R_0)_s$  it can be seen, by
using Lemma \ref{tbasic}(ii), that the relation
$(x^-_{\alpha_i}\otimes t^2)\vsim=0$ is actually a consequence of
$(y_{\alpha_i}^-\otimes t)v^\sigma_{i,m}=0.$
The modules $KR^\sigma(m\omega_i)$ are clearly graded
 modules for the graded algebra $\lie g[t]^\sigma$. Any $\lie g_0$--module  $V$ can be regarded as a module for $\lie g[t]^\sigma$ by setting
 $$(x_\epsilon\otimes t^{2r+\epsilon})v=0,\ \ v\in V,\ \  x_\epsilon\in\lie g_\epsilon, r\in\bz_+,$$ and we denote the corresponding module by $\ev_0^\sigma(V)$. The next lemma is easily checked.

\begin{lem}\label{eval}  For all $i\in I_0$ and $m\in\bz_+$ the module ${\rm{ev}}_0^\sigma(V(m\omega_i))$ is a quotient of $\krsm$. \hfill\qedsymbol
\end{lem}

 \subsection{The graded character of $\krsm$}
For $i\in I_0$ and  $m\in\bn$, let $P^+_0(i,m)^\sigma$ be the subset of $P_0^+$ defined by
$$P^+_0(i,m)^\sigma=P^+_0(i,d_i^\sigma)^\sigma+ P^+_0(i,m-d_i^\sigma)^\sigma ,$$
where
\begin{alignat*}{3}
&d_i^\sigma =1,&&\ \ \lie g\ne A_{2n},\\
&d_i^\sigma = 2,&&\ \ i\ne n, \ \ d_n^\sigma =4,\ \ \lie g=A_{2n},
\end{alignat*}
and,  for $1\le m\le d_i^\sigma$, $P_0^+(i,m)^\sigma$
is given as follows. If $\lie g$ is of type $A_{2n}$, then   \begin{eqnarray*}
 P_0^+(i,1)^\sigma &=&\{\omega_i\}, \ \
\\ P^+_0(i,2)^\sigma&=&\{2\omega_i, 2\omega_{i-1},
\cdots, 2\omega_1,0\},\ \ \ i\ne n,\\
P^+_0(n,2)^\sigma&=&\{2\omega_n\},\
\\ P^+_0(n,3)^\sigma&=&\{3\omega_n\}, \\
P^+_0(n,4)^\sigma& =& \{4\omega_n,2\omega_{n-1},\cdots ,2\omega_1, 0\}.
 \end{eqnarray*}  If $\lie g$ is  of type $A_{2n-1}$, then,
$$P^+_0(i,1)^\sigma=\{\omega_i,\omega_{i-2},\cdots ,\omega_{\overline i}\},$$ where $\overline i\in\{0,1\}$ and $i=\overline i\mod 2$. Finally if $\lie g$ is of type $D_{n+1}$, then
\begin{align*}
P^+_0(i,1)^\sigma&=\{\omega_i,\omega_{i-1},\cdots, 0\},\ \ \ i\ne n,\\
P^+_0(n,1)^\sigma&=\{\omega_n\}.
\end{align*}
  In particular  $P^+_0(i,m)^\sigma$ is finite.
 Fix an enumeration $\mu_s$, $0\le s\le k$, of the sets $P^+_0(i,m)^\sigma$, $1\le m\le d_i^\sigma$
by requiring $$\mu_s-\mu_{s+1}\in R_1^+,\ \ \mu_s-\mu_{s+2}\notin (R_0^+\cup R_1^+).$$
For $\mu\in P^+_0(i,m)^\sigma$ define a reduced expression and $|\mu|\in\bz_+$ analogously to the untwisted  case.

\noindent The main result  is the following.
\begin{thm}\label{maintw}Let $i\in I_0$, $m\in\bz_+$. \begin{enumerit}
\item
 For all $s\in\bz_+$,
$$KR^\sigma(m\omega_i)[s]\cong_{\lie g_0} \opl_{  \mu\in P^+_0(i,m)^\sigma:|\mu|=s}^{}V(\mu).$$
\item Write $m=d_i^\sigma m_0+m_1$, where $0\le
m_1<d^\sigma_i$. The canonical homomorphism of $\lie g^\sigma[t]$--modules  $$
KR^\sigma(m\omega_i)\to KR^\sigma(m_1\omega_i)\otimes KR^\sigma(d_i^\sigma\omega_i)^{\otimes
m_0}$$ is injective.
\end{enumerit}
\end{thm}
The proof of the theorem proceeds  as in the untwisted case.

\subsection{Elementary Properties of $\krsm$.}\hfill

 The next proposition is the twisted version of Proposition
\ref{eleprop} and is proved in the same way.
  \begin{prop}\hfill\label{twprop}
 \begin{enumerit}\iffalse
 \item $KR^\sigma(0)\cong \bc$.
\item For $\alpha\in R^+_0$ and $\beta\in R^+_1$, we have
$$(x^-_\alpha\otimes
t^2)^{m\omega_i(h_\alpha)+1}v^\sigma_{i,m}=(y^-_\beta\otimes
t)^{m\omega_i(h_\beta)+1}v^\sigma_{i,m} = 0.$$\fi \item We have
$$KR^\sigma(m\omega_i)=\opl_{\mu\in\lie h_0^*}^{}
KR^\sigma(m\omega_i)_\mu $$ and  $KR^\sigma(m\omega_i)_\mu\ne
0$ only if $\mu\in m\omega_i-Q_0^+$.
\item Regarded as a $\lie g_0$--module, $KR^\sigma(m\omega_i)$ and
$KR^\sigma(m\omega_i)[s]$, $s\in\bz_+$, are isomorphic to a direct
sum of irreducible finite--dimensional representations of $\lie
g_0$. In particular if $W_0$ is the Weyl group of $\lie g_0$, then
$KR^\sigma(m\omega_i)_\mu\ne 0$ iff
$KR^\sigma(m\omega_i)_{w\mu}\ne 0$ for all $w\in W_0$.
\item Given $0\le r\le m$, there exists a canonical homomorphism
$KR^\sigma(m\omega_i)\to
KR^\sigma(r\omega_i)\otimes KR((m-r)\omega_i)$ of
graded $\lie g[t]^\sigma$--modules such that
$v^\sigma_{i,m}\mapsto v^\sigma_{i,r}\otimes v^\sigma_{i,m-r}$.
\end{enumerit}\hfill\qedsymbol
\end{prop}

\begin{cor} We have $$\krsm\cong \opl_{\mu\in P_0^+}^{}
V(\mu)^{m_\mu(i,m)}.$$ In particular $KR^\sigma(0)\cong\bc$.
\hfill\qedsymbol\end{cor}

\subsection{An upper bound for $m_\mu(i,m)$}

\begin{prop}\label{twupb} For all $i\in I_0$, $m\in\bz_+$, we have
$$\krsm\cong_{\lie g_0}\opl_{\mu\in P_0^+(i,m)^{\sigma}}^{} V(\mu)^{m_\mu},$$ where $m_\mu\in\{0,1\}$.
\end{prop}
\begin{cor} For all $i\in I_0$ and $m\in\bz_+$, the modules $\krsm$ are finite--dimensional.\hfill\qedsymbol\end{cor}
We postpone the proof of this proposition to the next section.

\subsection{An explicit construction of the modules $KR^\sigma(m\omega_i)$, $i\in I_0$,  $1\le m\le  d^\sigma_i$.}\hfill

The next lemma is standard and can be found in
\cite{FH}, \cite{V}. Let $\iota:\lie
g_0\to\wedge^2(\lie g_1)$ be the $\lie g_0$--module map such that $[\ ,\ ]\cdot\iota={\rm{id}}$.
 \begin{lem} \label{M}
 \begin{enumerit}
\item[(a)]  If $\lie g$ is of type $D_{n+1}$, $n\ge 2$, then  $\wedge^2(\lie g_1)\cong_{\lie g_0}\iota(\lie g_0)$.\item[(b)]   If $\lie g$ is of type $A_n$, $n\ne 3$, we have $\wedge^2(\lie g_1)\cong\iota(\lie g)\oplus V(\nu)$ , where \begin{enumerit}
\item $\nu=\omega_1+\omega_3$ if $\lie g$ is of type
$A_{2n-1}$, $n\ge 3$,
\item $\nu=2\omega_1+\omega_2$ if $\lie g$ is of type
$A_{2n}$, $n\ge 3$,
\item $\nu= 2\omega_1+2\omega_2$ if $\lie g$ is of type
$A_{4}$,
\item $\nu=6\omega_1$ if $\lie g$ is of type
$A_{2}$.
\end{enumerit}\end{enumerit}\hfill\qedsymbol\end{lem}

Assume that $V_s$,  $0\le s\le k$, are $\lie g_0$--modules, let $p_s:\lie g_1\otimes
V_s\to V_{s+1}$, $0\le s\le k-1$, be  $\lie g_0$--module
maps, and set $p_k=0$.  Set also $q_s=p_{s+1}(1\otimes p_s)\in\Hom_{\lie g_0}(\lie
g_1\otimes\lie g_1 \otimes V_s, V_{s+2}).$
Suppose that   one of the following two conditions hold:
\begin{enumerit}
\item[(a)]   \begin{equation}\label{conda} q_s(\wedge^2(\lie g_1)\otimes V_s)=0,\ \ \forall\ \ 0\le s\le  k-2,\end{equation}
\item[(b)]   for all $0\le s\le k-2$, $y,z,w\in
\lie g_1$, $v\in V_s$, we have
\begin{equation} \label{condb}q_s(V(\nu)\otimes V_s)=0,\ \ p_{s+2}(1\otimes
p_{s+1})(1\otimes p_s)
\left(((z\wedge w)\otimes y- y\otimes
(z\wedge w ))\otimes v\right)=0.\end{equation}
\end{enumerit}
\vskip 12pt

Let $x\in\lie g_0$, $y\in\lie g_1$, $v\in V_s$. The following
formulas define a graded  $\lie g[t]^\sigma$--module structure on
$V=\oplus_{s=1}^k V_s$. $$(x\otimes t^2)v =q_s(\iota(x)\otimes
v),\ \ \ \ \ (y\otimes t) v=p_s(y\otimes v), $$ and $$ (x\otimes
t^{2r})v = (y\otimes t^{2r-1})v=0,$$ for all $r\ge 2$ and $0\le
s\le k$. Furthermore,
$$ V[s]\cong_{\lie g_0}V_s.$$
Moreover, if  the maps $p_s$, $0\le s\le k-1$  are all surjective and $V_0=\bu(\lie g[t]^\sigma)v_0$, then  the resulting $\lie g[t]^\sigma$--module is cyclic on $v_0$.

 \begin{prop}\label{fundtw}  Let $i\in I_0$, $1\le m\le d^\sigma_i$. The modules
 $V(\mu_s)$, $0\le s\le k$, satisfy \eqref{conda} (resp. \eqref{condb}) if
 $\lie g$ is of type
 $A_n$ (resp. $D_n$). The resulting $\lie g[t]^\sigma$-- module is isomorphic to $\krsm$ and  $$\krsm[s]\cong_{\lie g_0}V(\mu_s).$$
\end{prop}

 \begin{pf}
 If $k=0$ there is nothing to prove
since it follows from Proposition \ref{twupb} and Lemma \ref{eval} that $KR^\sigma(m\omega_i)\cong\ev^\sigma_0(V(m\omega_i))$.
Assume that $k>0$.  Using Proposition \ref{tpgen} it is easy to see that
  $\Hom_{\lie
 g_0}(V(\mu_s)\otimes  V(\mu_{s+1}),\lie g_1)\ne 0$ and hence
 $$ \Hom_{\lie
 g_0}(\lie g_1\otimes V(\mu_s), V(\mu_{s+1}))\ne 0.$$

Suppose first  that $\lie g$ is of type $A_n$. Equation \eqref{conda} holds if we  prove that  $$\Hom_{\lie g_0} (\wedge^2(\lie g_1)\otimes V_s,V_{s+2}) =
 0,\ \ \forall \ \ 0\le s\le k-2.$$ Since $\mu_s-\mu_{s+2}\notin R_0$, it is immediate
that
 $\Hom_{\lie g_0}(\lie g_0\otimes V(\mu_s), V(\mu_{s+2}))=0$.
To prove that $\Hom_{\lie g_0}(V(\nu)\otimes V(\mu_s),
V(\mu_{s+2}))=0$, it suffices to prove that   $$\Hom_{\lie
g_0}(V(\mu_s), V(\nu)\otimes V(\mu_{s+2}))=0.$$   This is done
exactly as in the untwisted case by noting that
$\mu_{s+2}+\nu-\mu_s=\sum_{j=1}^{s-1}k_j\alpha_j$ for some
$k_j\in\bz_+$.

If  $\lie g$ is of type $D_{n+1}$, $n\ge 2$, then Lemma \ref{M} implies that the first condition in \eqref{condb} is trivially  satisfied. If $n=3
 $  the second condition  is also trivially
true since $k\le 3$.
If  $n>3$, let $N$ be the $\lie g_0$--submodule of  $T^3(\lie g_1)$
spanned by elements of the form $(x\otimes y-y\otimes x)\otimes
z-z\otimes(x\otimes y-y\otimes x)$, with $x,y,z\in\lie g_1$. Note that
 $$N\cap
\wedge^3(\lie g_1)=0,\ \ \ N\subset \left(\lie  g_1\otimes\wedge^2(\lie g_1)+\wedge^2(\lie g_1)\otimes\lie g_1\right).$$
 Now, it is not hard to see that,
 $$\lie
g_1\otimes\wedge^2(\lie
 g_1)\cong_{\lie g_0} V(\omega_1+\omega_2)\oplus V(\omega_3)\oplus V(\omega_1),\
  \ \ n>3,$$ and that
the $\lie g_0$--submodule $V(\omega_3)$ occurs in
$T^3(\lie g_1)$ with mulitplicity one. It  follows that  $$
N\cong V(\omega_1+\omega_2)^{\oplus r_1}\oplus
 V(\omega_1)^{\oplus r_2}, \ \ 0\le r_1, r_2\le 2.$$

We now prove that  $$\Hom_{\lie g_0}(N\otimes V(\mu_s), V(\mu_{s+3}))=0,\ 0\le s\le k-3,$$ which establishes the second condition in \eqref{condb}.
Since $(V(\mu_s)\otimes V(\mu_{s+3}))_{\omega_1}=0$, it follows that $\Hom_{\lie
g_0}(V(\omega_1)\otimes V(\mu_s), V(\mu_{s+3}))=0$ and we are left to show that
$$\Hom_{\lie g_0}(V(\omega_1+\omega_2)\otimes V(\mu_s),
V(\mu_{s+3}))=0.$$
For this it suffices to prove that $$\Hom_{\lie
g_0}(V(\mu_s), V(\omega_1+\omega_2)\otimes V(\mu_{s+3}))=0.$$ This
is done as usual by noting that
$\mu_{s+3}+\omega_i+\omega_2-\mu_s$ is in the span of the elements
$\alpha_r$, $1\le r\le s-1$
and by observing that the $\lie (g_0)_J$--module $\bu(\lie (g_0)_J)v_{\omega_1+\omega_2}\subset V(\omega_1+\omega_2)$ is not  dual to the $(\lie g_0)_J$--module  $\bu((\lie g_0)_J)v_{\mu_{s+3}}\subset V(\mu_{s+3})$, where $J=\{1,\cdots ,s-1\}$.

This proves that if we fix non--zero maps $p_s\in\Hom_{\lie g_0}(\lie g_1\otimes V(\mu_s), V(\mu_{s+1}))$, then we can construct a graded  cyclic $\lie g[t]^\sigma$--module
$$V=\bu(\lie g[t])v_{\mu_0}=\oplus_{s=0}^k  V(\mu_s),\ \ V[s]\cong V(\mu_s),\ \ 0\le s\le k.$$
As in the untwisted case, to complete the proof it suffices, in view of Proposition \ref{twupb}, to prove that
$V$ is a quotient of $\krsm$. Since
$\mu_0=m\omega_i$ and  $\mu_0-\mu_1\in R_0^+\backslash\{\alpha_i\}$
we get $$p_0\left((\lie n_1^+\oplus \lie h_1\oplus\bc
y^-_{\alpha_i})\otimes v_{m\omega_i}\right)=0,$$   i.e.,
$$(\lie n^+_1\otimes t^{2r-1}) v_{m\omega_i}=(\lie h_1\otimes
t^{2r-1})v_{m\omega_i}= (y^-_{\alpha_i}\otimes t)v_{m\omega_i}=0,\ \
\ r\ge 1. $$ For the same reasons, $$q_0\left((\lie n^+_0\oplus\lie h_0\oplus
\bc x^-_{\alpha_i})\otimes v_{m\omega_i}\right)=0.$$  Hence we get
$$(\lie n_0^+\otimes t^{2r})v_{m\omega_i}=(\lie h_0\otimes
t^{2r})v_{m\omega_i}=(x^-_{\alpha_i}\otimes t^{2r})\vsim=0,\ \ r\ge
1.$$ The relations \eqref{deftw2} follow since $V$ is obviously finite--dimensional
  which completes the  proof  that $V$ is a quotient of $\krsm$.
\end{pf}

Given $\mu\in P_0^+(i,m)^\sigma$, $m= d_i^\sigma m_0+m_1$, $0\le m_1<  d_i^\sigma$,  with reduced expression $$\mu=m_1\omega_i+\mu_{j_1}+\cdots +\mu_{j_{m_0}},$$ set $s_j=\#\{r:j_r\ge j\}$, $0\le j\le k$, and
 $$\boy_\mu=(y^-_{\mu_{k-1}-\mu_k}\otimes t)^{s_k}\cdots (y^-_{\mu_0-\mu_1}\otimes t)^{s_1}\in\bu(\lie n^-[t]^\sigma).$$
It is easily  seen that  $$\boy_\mu v_{i,m}^\sigma\in KR(m\omega_i)[|\mu|]\cap KR(m\omega_i)_{\mu}.$$
The next corollary is proved in the same way as Corollary \ref{fundun} and we omit the details.
\begin{cor} \label{fundtw}

For  $0\le s\le k$ we have
$$\boy_{\mu_s} v^\sigma_{i, d_i^\sigma}\ne 0,\ \ h(\boy_{\mu_s}
v^\sigma_{i, d_i^\sigma})=\mu_s(h) (\boy_{\mu_s}v^\sigma_{i, d_i^\sigma}),\ \ \lie
n^+(\boy_{\mu} v^\sigma_{i,d_i^\sigma})=0, \ \ h\in\lie h_0,$$ or, equivalently, $$\boy_{\mu_s}v_{\mu_0}=v_{\mu_s}.$$
\hfill\qedsymbol
\end{cor}

\subsection{The modules $\tilde{KR}^\sigma(m\omega_i)$ and the  completion of the proof of Theorem \ref{maintw}}\hfill

 For $m=d_i^\sigma m_0+m_1\in\bz_+$, $i\in I_0$,
 set
$$\tilde{KR}^\sigma(m\omega_i)=\bu(\lie g[t]^\sigma)(
v^\sigma_{m_1\omega_i}\otimes (v^\sigma_{d_i^\sigma\omega_i})^{\otimes
m_0})\subset {KR}^\sigma(m_1\omega_i) \otimes
{KR}^\sigma(d_i^\sigma\omega_i)^{\otimes m_0}.$$
The next  proposition is proved in exactly the same way as the corresponding result in Section 2.7 for
the untwisted case, by using Proposition \ref{fundtw} and Corollary \ref{fundtw},
and clearly completes the proof of Theorem \ref{maintw}.
\begin{prop}
 Let $i\in I_0$, $m\in\bn$. Then $$\tilde{KR}^\sigma(m\omega_i)[s]\cong_{\lie g_0}\bigoplus_{\mu\in P^+_0(i,m)^\sigma:|\mu|=s}V(\mu).$$
In particular,
$$\tilde{KR}^\sigma(m\omega_i)\cong KR^\sigma(m\omega_i)$$ as $\lie g[t]^\sigma$--modules.\hfill\qedsymbol

\end{prop}

\section{ Proof of Proposition \ref{twupb}}
We will use the following remark repeatedly without further
comment. Let $x\in(\lie g_\epsilon)_\alpha$, $\alpha\in R_\epsilon$, and $r\in\bz_+$ be such that $(x\otimes t^{2r+\epsilon})\vsim=0$. Then $(x\otimes t^{2r+2s+\epsilon})v_{i,m}^\sigma=0$
for all $s\in\bz_+$. To see this, observe that if $h\in\lie h_\epsilon$, then by definition we have $(h\otimes t
^{2s+\epsilon})v^\sigma_{i,m}=0$. If $x\in(\lie
g_\epsilon)_\alpha$ for some $\alpha\in R_\epsilon$,  choose
$h\in\lie h_0$ such that $[h,x]= x$. This gives $$ 0= [h\otimes
t^2,x\otimes t^r]v^\sigma_{i,m} = (x\otimes
t^{r+2})v^\sigma_{i,m},$$ thus proving the remark.
\subsection{The case $\alpha_i\in R_1$ with $\varepsilon_i(\phi)=1$} The following Lemma establishes the proposition in this case.
\begin{lem}\label{indbegin} Suppose that $i\in I_0$ is such that $\alpha_i\in R_1$ and $\varepsilon_i(\phi)=1$ (in other words $i=n$ if $\lie g$ is of type $D_{n+1}$ and $i=1$ if $\lie g$ is of type $A_{2n-1}$).  Then $$\krsm\cong\ev_0^\sigma(V(m\omega_i)).$$\end{lem}
\begin{pf} It is straightforward to check that   we can write \begin{equation}\label{br}y_\phi^-=[x^-_\beta, y^-_{\alpha_i}],\ \ x^-_\theta=[y^-_\phi,y^-_\gamma],\end{equation} for some $\beta\in R_0^+$, $\gamma\in R_1^+$  with $\varepsilon_i(\beta)=0$, $\varepsilon_i(\gamma)=1$.
 It follows from \eqref{deftw4} and Proposition \ref{twprop} that $(y_\phi^-\otimes t)\vsim=0.$ Since for all $\alpha\in R_1^+$ we have  $y_\alpha^-\otimes t\in \bu(\lie n^+)(y^-_\phi\otimes t)$ it follows from \eqref{deftw1} that $(y^-_\alpha\otimes t)\vsim =0$. It follows now from \eqref{br} that $(x^-_\theta\otimes t^2)\vsim=0$ and hence that $(x^-_\alpha\otimes t^2)\vsim=0$ for all $\alpha\in R_0^+$. This proves that $\krsm=\bu(
 \lie n_0^-)\vsim$ and the Lemma is proved.

\end{pf}
\subsection{The subalgebras $\lie n^-[t]^\sigma_{\max(i)}$}\hfill

From now on we assume that $i\in I_0$ is such that either $\alpha_i\notin R_1$ or $\varepsilon_i(\phi)=2$.
Let $\max(i)\in\{1,2\}$ be the  maximum value of the restriction of $\varepsilon_i$ to $ R_1^+$, and let
$\lie n^-[t]^\sigma_{\max(i)}$ be the subalgbera of $\lie n^-[t]^\sigma$ generated by $$\{y_\alpha^-\otimes t:\alpha
\in R^+_1,\ \ \varepsilon_i(\alpha)=\max(i)\}.$$
If $\lie g$ is of type $A_n$, then it is easy to see that
\begin{equation}\label{an1} \lie n^-[t]^\sigma_{\max(i)}=\bigoplus_{\alpha\in R_1:\varepsilon_i(\alpha)=\max(i)}\bc( y^-_\alpha \otimes t),\end{equation}
while if $\lie g$ is of type $D_{n+1}$, we have

\begin{equation}\label{dn1i} \lie n^-[t]^\sigma_{\max(i)}
=\left(\opl_{\alpha\in R_1:\varepsilon_i(\alpha)=1}^{}\bc (y^-_\alpha\otimes t)\right)\oplus\left( \opl_{\alpha\in R_0^+:\varepsilon_i(\alpha)=2}^{}\bc (x^-_\alpha\otimes t^2)\right).\end{equation}

\subsection{Further relations satisfied by $\vsim$}\hfill

\begin{prop}\label{span1}\hfill
\begin{enumerit}
\item Let $\alpha\in R_1^+$. If $\varepsilon_i(\alpha)<\max(i)$ we have
 $$(y_\alpha^-\otimes t^{2r+1})\vsim=0 \qquad\text{for all}\qquad r\in\bz_+,$$
while if $\varepsilon_i(\alpha)=\max(i)$ we have
 $$(y_\alpha^-\otimes t^{2r+3})\vsim=0 \qquad\text{for all}\qquad r\in\bz_+.$$
\item If $\lie g$ is of type $A_n$, then for all $\alpha\in R_0^+$, $r\in\bz_+$, we have $$(x^-_\alpha\otimes t^{2r+2})\vsim=0.$$
\item If $\lie g$ is of type $D_{n+1}$, then for all $\alpha\in R_0^+$ and $r\in\bz_+$, we have $$(x^-_\alpha\otimes t^{2r+2\varepsilon_i(\alpha)})\vsim=0,\ \ i<n.$$
\end{enumerit}
\end{prop}
\begin{pf}Let $\alpha\in R_1^+$.  Suppose that $\varepsilon_i(\alpha)=0$. Then
 $s_\alpha(m\omega_i-\alpha)=m\omega_i+\alpha$ and hence $(\krsm)_{m\omega_i-\alpha}=0$.
 This shows that $$x^-_\alpha\vsim =(y_\alpha^-\otimes t)\vsim=0.$$
  If $\varepsilon_i(\alpha)=1<\max(i)$, then  we can choose $j\in I_0$ such that
 $\beta=\alpha-\alpha_j\in R_1^+\cap R_0^+$. If $j=i$, then\  $\varepsilon_i(\beta)=0$ and we get by using \eqref{deftw4} that $$(y^-_\alpha\otimes t)\vsim=[x^-_\beta,y^-_{\alpha_i}\otimes t]\vsim=0,$$ and if $j\ne i$, then $\varepsilon_i(\alpha_j)=0$ and we get $$(y^-_\alpha\otimes t)\vsim=[x^-_{\alpha_j},y^-_{\beta}\otimes t]\vsim = x^-_{\alpha_j}(y^-_{\beta}\otimes t)\vsim.$$
 Repeating the argument with $\beta$ we eventually get $(y^-_\alpha\otimes t)\vsim=0,$ thus proving (i). The proofs of (ii) and (iii) are similar and we omit the details.
\end{pf}

\begin{cor}\label{spanspeccase}\hfill
\begin{enumerit}
\item As vector spaces we have
\begin{equation}\label{span2}\krsm=\bu(\lie n^{-}_0)\bu(\lie n^-[t]^\sigma_{\max(i)})\vsim.
\end{equation}  In particular $\krsm$ is finite--dimensional.
\item  If   $\lie g$ is of type $A_2$ and $1\le m\le 3$,
 then we have $$\krsm\cong{\rm{ev}}_0^\sigma(V(m\omega_i))$$ and Proposition \ref{twupb} is proved in these cases.\end{enumerit}
\end{cor}
\begin{pf} A straightforward application of the Poincar\'e-Birkhoff-Witt theorem gives \eqref{span2}.
Since   $\lie n^{-}_0\oplus  \lie n^-[t]^\sigma_{\max(i)}$ is a finite--dimensional Lie algebra,
it follows that for all $\mu\in P$ $$\dim(\krsm_\mu)<\infty.$$  Proposition \ref{twprop} now implies  that
 if $\krsm$ is infinite--dimensional, there must be infinitely many elements $\nu_r\in P^+$ with $\krsm_\nu\ne 0$ and hence $\nu_k\in P^+\cap (m\omega_i-Q^+)$. But this is a contradiction since it is well-known that the set $P^+\cap (m\omega_i-Q^+)$ is finite.

  If $\lie g$ is of type $A_2$ and $m\le 3$, then $i=1$ and   $(m\omega_1-\phi)(h_{\alpha_1})<0$ since $\varepsilon_1(\phi)=2$.
On the other hand, since $\varepsilon_1(\phi-\alpha_1)=1$, we see that
 $$x^+_{\alpha_1}(y^-_\phi\otimes t)\vsim= (y^-_{\phi-\alpha_1}\otimes t)\vsim =0.$$ Since $\krsm$ is isomorphic to a direct  sum of finite--dimensional $\lie g_0$--modules this forces $(y^-_\phi\otimes t)\vsim=0$ which proves the corollary.
\end{pf}
\subsection{Proof of Proposition \ref{twupb}}
Recall the enumeration $\mu_0,\cdots ,\mu_k$ of the sets $P^+_0(i, d_i^\sigma)^\sigma$ and set $\phi_s=\mu_{k-s}-\mu_{k-s+1}\in R_1^+$ for $1\le s\le k$. Given $\bor=(r_1,\cdots , r_k)\in\bz_+^k$, set
 $\eta_\bor=\sum_{s=1}^kr_s\phi_s\in Q^+$  and $$\boy_\bor=(y_{\phi_1}\otimes t)^{r_1}\cdots (y_{\phi_k}\otimes t)^{r_k}.$$
Note that $\eta_\bor+\eta_{\bor'}=\eta_{\bor+\bor'}$ for all $\bor, \bor'\in \bz_+^k$. For $1\le s\le k$, let  $\boe_s\in\bz_+^k$ be the element with one in the $s^{th}$ place and zero elsewhere.
As usual $\le$ is the partial order on $P_0$ defined by $\mu\le \mu'$ iff $\mu'-\mu\in Q_0^+$.
\begin{prop} \label{mainredtry} We have,  \begin{equation}\label{hw} \lie n ^+_0\boy_\bor\vsim\in\sum_{\eta_\bor'<\eta_\bor}\bu(\lie n^-_0)\boy_{\bor'}\vsim.\end{equation}
\begin{equation}\label{gspan}\krsm=\sum_{\bor\in\bz_+^k}\bu(\lie g_0)\boy_\bor\vsim.\end{equation}

\end{prop}
\nc\pr{{\text{pr}}}

Assuming this proposition we complete the proof of Proposition \ref{twupb}  as follows.
By Proposition \ref{twprop} we can pick  a $\lie g_0$--module $W^0$ such that $$\krsm\cong_{\lie g_0}\bu(\lie g_0)\vsim \oplus W^0.$$ If $W^0=0$ there is nothing to prove.  Otherwise, let $\pr_{W^0}$ be the $\lie g_0$-projection
of $\krsm$ onto $W^0$. Using  \eqref{gspan} we  see that there exists $\bor\in\bz_+^k$ such that $\pr_{W^0}(\boy_\bor)\ne 0$.
Choose $\bor_1\in\bz_+^k$ such that $v_{\bor_1}=\pr_{W^0}(\boy_{\bor_1}\vsim)\ne 0$ and such that
$\eta_{
\bor_1}$ is minimal: i.e., if $\pr_{W^0}(\boy_{\bor}\vsim)\ne 0$
for some $\bor\in\bz_+^k$, then
$\eta_{\bor_1}-\eta_{\bor}\notin Q_0^+$.
Using \eqref{hw}  we get $$\lie n^+_0\pr_{W^0}(\boy_{\bor_1}\vsim)=\pr_{W^0}(\lie n_0^+\boy_{\bor_1}\vsim) =0,$$ i.e., $$\lie n_0^+(v^\sigma_{\bor_1})=0.$$  Thus we can write $$\krsm\cong_{\lie g_0}\bu(\lie g_0)\vsim \oplus \bu(\lie g_0) v^\sigma_{\bor_1}\oplus W^1,$$ for some $\lie g_0$--submodule $W^1$.
Repeating the argument we find that there exists $\bor_s$, $s\in\bz_+$, such that $$\krsm\cong\bu(\lie g_0)\vsim \oplus\opl_{s\ge 1}^{}\bu(\lie g_0)v^\sigma_{\bor_s},$$ with $\lie n_0^+v^\sigma_{\bor_s}=0$.
 Since  $$v^\sigma_{\bor_s}\in\krsm_{m\omega_i-\sum_{j=1}^k s_j\phi_j},$$ for some $s_j\in\bz_+$ we get
$$m\omega_i-\sum_{j=1}^k s_j\phi_j=(m- s_kd^\sigma_i)\omega_i+ (s_k-s_{k-1})\mu_1+\cdots  +
(s_{2}-s_{1})\mu_{k-1}+s_1\mu_k\in P_0^+.$$  An inspection of the sets $P^+_0(i,d_i^\sigma)^\sigma$ shows that this implies  $$m\ge d_i^\sigma s_k,\ \ s_j\le s_{j+1},\ \ 1\le j\le k-1,$$ which in turn implies that the weight of $v^\sigma_{\bor_s}$ is in $P^+_0(i,m)^\sigma$. Since $\bor_s\ne\bor_{s'}$ if $s\ne s'$, Proposition \ref{twupb} is now proved.
It remains to prove Proposition \ref{mainredtry}.

\subsubsection{Proof of \eqref{hw}}
It is useful to recall  our convention that if $\alpha \in R_1^+$, then  $y_\alpha^-$ denotes an arbitrary non--zero element of the space $(\lie g_1)_{-\alpha}$.

{\bf{ The case when $\lie g$ is of type $D_{n+1}$ or $A_{2n}$}.}\hfill

We proceed    by induction on $n$. To see that  induction begins for the series $D_{n+1}$ at $n=2$ note that then,  $\lie n^+[t]^\sigma_{\max(1)}=\bc(y^-_{\alpha_1+\alpha_2}\otimes t)$ and hence the result follows from Corollary \ref{spanspeccase}(i).  For  $A_{2n}$ it begins
  at $n=1$ by noting that
 $\lie n^-_0[t]^\sigma_{\max(1)}=\bc (y_{\phi}\otimes t)$ and $\varepsilon_1(\phi-\alpha_1)=1<\max(1)$ and again using Corollary \ref{spanspeccase}.

For the inductive step set $J=\{2,\cdots ,n\}$ and let $j\in J$.
We have
$$x^+_j\boy_\bor\vsim=(y_{\phi_1}\otimes t)^{r_1}x_j^+\boy_{\bor-r_1\boe_1}\vsim\in (y_{\phi_1}\otimes t)^{r_1}\sum_{\bor'\in S}\bu(
(\lie n^-_0)_J)\boy_{\bor'}\subset \sum_{\bor'\in S}\bu(
(\lie n^-_0)_J)(y_{\phi_1}\otimes t)^{r_1}\boy_{\bor'}
$$ where $$S=\{\bor'\in\bz_k^+:r_1'=0, \eta_{\bor'}<
\eta_{\bor-r_1\boe_1}\}.$$ The last inclusion is a consequence of the fact
 that $[y_{\phi_1}^-,\lie n^-_0]=0$. It is now simple to check that $
\eta_
\bor> \eta_{\bor'+r_1\boe_1}$.

Suppose now that $j=1$ and that $\lie g$ is of type $A_{2n}$.
We have
$$x^+_{\alpha_1}\boy_\bor\vsim =(y^-_{\alpha_1+2\sum_{s=2}^n\alpha_s}\otimes t)\boy_{\bor-\boe_1}\vsim= \delta_{i,1} x^-_{\alpha_1}\boy_{\bor-\boe_1+\boe_2}\vsim,$$
where $\delta_{i,1}$ is one if $i=1$ and zero otherwise. If  $\lie g$ is of type $D_{n+1}$,
we have \begin{equation}\label{step1} x_1^+\boy_\bor\vsim=\boy_{\bor-\boe_1+\boe_2}\vsim+(x^-_\theta\otimes t^2)\boy_{\bor-2\boe_1}\vsim,\end{equation} where we recall that
$\theta=\alpha_1+2\sum_{j=2}^n\alpha_j$.
Since  $\eta_{\bor-\boe_1+\boe_2}<\eta_{
\bor}$  it remains to prove that $$(x^-_\theta\otimes t^2)\boy_{\bor-2\boe_1}\vsim\in\sum_{\eta_{\bor'}<\eta_{\bor}}\bu(\lie n^-_0)\boy_{\bor'}\vsim.$$
Now\begin{equation}\label{fstep} x^-_{\alpha_1}\boy_{\bor-2\boe_1+2\boe_2}\vsim= (x^-_{\theta}\otimes t^2)\boy_{\bor-2\boe_1}\vsim +\boy_{\bor-\boe_1+\boe_2}\vsim. \end{equation}  Since $\eta_{\bor}>\eta_{\bor-\boe_1+\boe_2}$  and $\eta_{\bor}>\eta_{\bor-2\boe_1+2\boe_2}$  the inductive step is proved.

\nc\End{{\text{End}}}

\vskip 18pt

{\bf{The case when $\lie g$ is of type $A_{2n-1}$}.}

 We note that  induction begins at $n=2$ by using Corollary \ref{spanspeccase}, since in this case we have when $i=2$ that $\lie n^+[t]^\sigma_{\max(2)}=\bc(y^-_{\alpha_1+\alpha_2}\otimes t)$.
Set $J=\{2,,\cdots ,n\}$. For the inductive step, note that if $i$ is odd, then $\phi_s\in R_J^+$ and by the induction hyopthesis it suffices to consider only the case $x^+_{\alpha_1}\boy_\bor\vsim$. Since   $$[x^+_{\alpha_1}, \boy_\bor]=0,$$ the result is immediate.

Assume then that $i$ is even and let $j\ne 2$. We have
$$x^+_j\boy_\bor\vsim=(y_{\phi_1}\otimes t)^{r_1}x_j^+\boy_{\bor-r_1\boe_1}\vsim\in (y_{\phi_1}\otimes t)^{r_1}\sum_{\bor'\in S}\bu(
(\lie n^-_0)_J)\boy_{\bor'}\subset \sum_{\bor'\in S}\bu(
(\lie n^-_0)_J)(y_{\phi_1}\otimes t)^{r_1}\boy_{\bor'}
$$ where $$S=\{\bor'\in\bz_k^+:r_1'=0, \eta_{\bor'}<
\eta_{\bor-r_1\boe_1}\}.$$
If $j=2$, we get  $$x^+_{\alpha_2}\boy_\bor\vsim =(y^-_{\alpha_1+\alpha_2+2\sum_{s=3}^n\alpha_s}\otimes t)\boy_{\bor-\boe_1}\vsim= \delta_{i,2} x^-_{\alpha_1+\alpha_2+\alpha_3}\boy_{\bor-\boe_1+\boe_2}\vsim.$$ This completes the proof of the inductive step.

\subsubsection{Proof of \eqref{gspan}} Let $\lie k_i^+\subset \lie n^+_0$ be the subalgebra spanned by elements $\{x^+_\alpha:\varepsilon_i(\alpha)=0\}$.
It is easily seen that the subalgebra $\lie n^-[t]^\sigma_{\max(i)}$ is a module for $\lie k_i$ via the adjoint action, i.e., $$[\lie k_i,\lie n^-[t]^\sigma_{\max(i)}]\subset \lie n^-[t]^\sigma_{\max(i)}.$$ This action  defines a $\lie k_i$--module structure on  $\bu(\lie n^-[t]^\sigma_{\max(i)} )$ and  let $\rho:\lie k_i\to\End(\bu(\lie n^-[t]^\sigma_{\max(i)}))$ denote  the corresponding homomorphism of Lie algebras.

\begin{lem}\label{mod}  We have $$\bu(\lie n^-[t]^\sigma_{\max(i)})=\sum_{\bor\in\bz_+^k}\rho(\bu(\lie k_i^+))\boy_\bor.$$
\end{lem}
Assuming the Lemma we complete the proof of \eqref{gspan}  as follows. Let $g\in\bu(\lie n^-[t]^\sigma_{\max(i)})$ and write $g=\sum_{\bor\in\bz_k^+}\rho(\bx_\bor)\boy_\bor$ for some $\bx_\bor\in\bu(\lie k_i^+)$.
Then,  $$ g\vsim=\sum_{\bor\in\bz_+^k}\rho(\bx_\bor)\boy_r\vsim=\sum_{\bor\in\bz_+^k}\bx_\bor\boy_\bor\vsim\in\sum_{\bor\in\bz_+^k}\bu(\lie g_0)\boy_r\vsim,$$ where the second equality uses the fact that $\bx_\bor\vsim=0$. An application of Corollary \ref{spanspeccase}(i) now  gives \eqref{gspan}.

The proof of Lemma \ref{mod} is  straightforward and we give a proof when $\lie g$ is of type $D_{n+1}$. The other cases are similar. If $i=1$, there is nothing to prove and we assume from now on that $1<i<n$.
We proceed by induction on $n$, with induction beginning at $n=2$ by the preceding comment. Let $J=\{2,\cdots ,n\}$, by the induction hypothesis we have
$$\bu(\lie n^-[t]^\sigma_{\max(i)\cap\lie n^-_J[t]^\sigma})=\sum_{\bor\in\bz_+^i:\  r_1=0}\rho(\bu(\lie k_i^+\cap\lie n^+_J))\boy_\bor,$$ and since $[y^-_{\phi_1},\lie n^-[t]^\sigma]= [y^-_{\phi_1}, \lie n^+_J]=0$ we get
\begin{equation}\label{reda}\sum_{r_1\in\bz_+}(y_{\phi_1}\otimes t)^{r_1}\bu(\lie n^-[t]^\sigma_{\max(i)\cap\lie n^-_J[t]^\sigma})=\sum_{\bor\in\bz_+^i}\rho(\bu(\lie k_i^+\cap\lie n^+_J))\boy_\bor.\end{equation}
Thus the lemma is established if we prove that
\begin{equation}\label{red}\rho(\bu(\lie k_i^+))\sum_{r_1\in\bz_+}(y_{\phi_1}\otimes t)^{r_1}\bu(\lie n^-[t]^\sigma_{\max(i)\cap\lie n^-_J[t]^\sigma})=\bu(\lie n^-[t]^\sigma_{\max(i)}).
\end{equation}
For $1\le j\le i-1$, set
$$\theta_{j} = \alpha_1+\cdots+\alpha_{j}+2(\alpha_{j+1}+\cdots+\alpha_n),$$ and  $$\bx_\bos=(x^-_{\theta_1}\otimes t^2)^{s_1}\cdots (x^-_{\theta_{i-1}}\otimes t^2)^{s_{i-1}},\ \ \bos\in\bz_+^{i-1}. $$ By the PBW theorem we have,
$$\bu(\lie n^-[t]^\sigma_{\max(i)})=\sum_{\bos\in\bz_+^{i-1}, \ r\in\bz_+}\bx_\bos (y^-_{\phi_1}\otimes t)^r\bu(\lie n^-[t]^\sigma_{\max(i)\cap\lie n^-_J[t]^\sigma}),$$
and hence \eqref{red} follows if we prove that for all $\bos\in\bz_+^{i-1}$ and $r\in\bz_+$ we have
\begin{equation}\label{redb}\bx_\bos (y^-_{\phi_1}\otimes t)^r\bu(\lie n^-[t]^\sigma_{\max(i)\cap\lie n^-_J[t]^\sigma})\in \rho(\bu(\lie k_i^+))\sum_{r_1\in\bz_+}(y_{\phi_1}\otimes t)^{r_1}\bu(\lie n^-[t]^\sigma_{\max(i)\cap\lie n^-_J[t]^\sigma}).\end{equation}

To prove \eqref{redb} let $r_1\in\bz_+$, $g\in \bu(\lie n^-[t]^\sigma_{\max(i)\cap\lie n^-_J[t]^\sigma})$. For $2\le 2s_1\le r_1$ we get
$$\rho(x^+_{\alpha_1})^{s_1} (y_{\phi_1}\otimes t)^{r_1}g= \sum_{p=0}^{s_1}(x^-_{\theta_1}\otimes t^2)^{p}(y_{\phi_1}\otimes t)^{r_1-s_1-p}(y_{\phi_2}\otimes t)^{s_1-p}g , $$ which proves \eqref{redb} for all $\bos$ with $s_j=0$ for $j>1$.

Assume that \eqref{redb} is established for all $\bos$ such that $s_j=0$ for all $j\ge r$. To prove that it holds for $\bos$ with $s_j=0$ for $j>r>1$,
suppose first that $s_r=1$. Then, $$\rho(x^+_{\alpha_r})\bx_{\bos -\boe_r+\boe_{r-1}}(y_{\phi_1}\otimes t)^{r_1}g= \bx_\bos (y_{\phi_1}\otimes t)^{r_1}g+   \bx_{\bos-\boe_r+\boe_{r-1}} (y_{\phi_1}\otimes t)^{r_1}\rho(x^+_{\alpha_r})g,$$
and hence we get by using \eqref{reda} that \eqref{redb} holds for $\bx_\bos (y_{\phi_1}\otimes t)^{r_1}g$ when $s_r=1$. An obvious induction on $s_r$ gives the result in general.

\section{Demazure modules and concluding remarks}
\subsection{The relation between the modules $KR(m\omega_i)$ and Demazure modules}

\subsubsection{The affine Kac--Moody algebras}
  \def \whg{{\widehat{\lie g}}}
 \def \whh{{\widehat{\lie h}}}

Let  $\bc[t,t^{-1}]$ be the ring of Laurent polynomials in an
indeterminate $t$. Let $\whg$ be the affine Lie algebra defined by
$$\widehat{\lie g}=\lie g\otimes\bc[t,t^{-1}]\oplus\bc c\oplus \bc
d,$$ where $c$ is central and $$[x\otimes t^s, y\otimes
t^k]=[x,y]\otimes t^{s+k}+s\delta_{s+k,0}\left<x,y\right>c,\ \
[d,x\otimes t^s]=sx\otimes t^s,$$ for all $x,y\in\lie g$,
$s,k\in\bz$ and where $\left<\ ,\ \right>$ is the Killing form of
$\lie g$. Set $$\whh=\lie h\oplus\bc c\oplus\bc d, $$ and regard
$\lie h^*$ as a subspace of  $\whh^*$ by setting
$\lambda(c)=\lambda(d)=0$  for $\lambda \in \lie h^*$.
For $0\le i,j \le n$ let $\Lambda_i\in\whh^*$ be defined by
$\Lambda_i(h_j)= \delta_{i,j},$  and let $\widehat{P}^+\subset
\whh^*$ be the non--negative integer linear span of $\Lambda_i$,
$0\le i\le n$. Define $\delta\in(\whh)^*$ by $$\delta|_{\lie h} =0, \ \ \delta(d)=1,\ \ \delta(c)=0.$$  The elements  $\alpha_i\in(\whh)^*$,
$0\le i\le n$, where $\alpha_0=\delta-\theta$ are  a set of simple roots for $\whg$ with respect to $\whh$.
Let $\widehat{Q}^+$ be the subset of $\widehat{P}^+$ spanned by
the elements $\alpha_i$, $0\le i\le n$. Let $\widehat{W}$ be the
extended affine Weyl group and $(\ ,\ )$ be the $\widehat W$--invariant form on $\whh
^*$ obtained by requiring $(\Lambda_i,\alpha_j)=\delta_{ij} $ for
all $0\le i,j\le n$ and $(\delta,\alpha_i)=0$ for all $0\le i\le
n$.
\subsubsection{The highest weight integrable modules}

Given $\Lambda=\sum_{i=0}^n m_i\Lambda_i\in\widehat{P}^+$, let
$L(\Lambda)$ be the $\whg$--module generated by an element
$v_\Lambda$ with relations: $$\lie gt[t] v_\Lambda=0,\ \quad (\lie
n^+ \otimes 1)v_\Lambda=0,\quad (h \otimes 1) v_\Lambda=\Lambda(h)
v_\Lambda, $$ $$(x_{\alpha_i}^- \otimes 1)^{m_i+1}v_\Lambda=0,\ \
(x_{\theta}^+\otimes t^{-1})^{m_0+1}v_\Lambda=0,$$ for $1\le i\le n$.
This module is known to be irreducible and integrable (see
\cite{K}). The next proposition also can be found in \cite{K}.
\begin{prop} \label{wtL} Let $\Lambda\in \widehat{P}^+$.
\begin{enumerit}
\item  We have
$$L(\Lambda)=\bigoplus_{\Lambda'\in
\widehat{P}}L(\Lambda)_{\Lambda'},\qquad \mbox{where} \quad
L(\Lambda)_{\Lambda'}=\{v\in L(\Lambda): hv=\Lambda'(h)v,\ \
h\in\whh\}.$$ Moreover, $\dim(L(\Lambda)_{\Lambda'})<\infty$,
$\dim(L(\Lambda))_{\Lambda}=1$, and  $L(\Lambda)_{\Lambda'}=0$ if
$\Lambda-\Lambda'\notin\widehat{Q}^+$.
\item The set
\begin{equation*} \wt(L(\Lambda))=\{\Lambda'\in \widehat{P}:
L(\Lambda)_{\Lambda'}\ne 0\}\subset \Lambda-\widehat{Q}^+.
\end{equation*}
is preserved by $\widehat{W}$ and $$\dim(L(\Lambda)_{w\Lambda'} )=
\dim(L(\Lambda)_{\Lambda'} )\ \ \forall\ \ w\in\widehat{W},\ \
\Lambda'\in \widehat{P}.$$ In particular,
$\dim(L(\Lambda)_{w\Lambda} )=1$ for all $w\in\widehat{W}$.
\end{enumerit}\hfill\qedsymbol
\end{prop}

\subsubsection{The Demazure modules } From now on we fix $\Lambda=m\Lambda_0$ for some $m\in\bz_+$ and we also assume for simplicity that $\lie g$ is of tyep $A_n$ or $D_n$.  Given $w\in\widehat{W}$, let $v_{w\Lambda}$ be a
non--zero element in $L(\Lambda)_{w\Lambda}$ and let
$$D(w\Lambda)=\bu(\lie g[t])v_{w\Lambda}\subset L(\Lambda).$$
 It is clear from Proposition~\ref{wtL} that $D(w\Lambda)$ is
finite--dimensional for all $w\in \widehat{W}$.
\begin{prop}\label{demquot}
Let $\Lambda=m\Lambda_0
\in\hat P^+$ for some $m>0$ and $w\in\widehat{W}$ be such that $w\Lambda|_{\lie h}=m\omega_i$ for some $1\le i\le n$. Then, $D(w\Lambda)$ is isomorphic to
$KR(m\omega_i)$.
\end{prop}
\begin{pf}
To prove the proposition we first   show that $v_{w\Lambda}$
satisfies the defining relations of $KR(m\omega_i)$.
It was shown in \cite{CL} that the element $v_{w\Lambda}$ satisfied the relations $$\lie n^+[t] v_{w\Lambda}=0,\ \ (h\otimes t^r)v_{w\Lambda}=0, \ \ h\in\lie h, r\in\bn.$$
The proof given in \cite[Proposition 1.4.3]{CL} was for the special case when $\lie g$ is of type $A_n$, but     works for any simple Lie algebra. Since $D(w\Lambda)$ is finite--dimensional the only relation that remains to be checked is that $$(x^-_{\alpha_i}\otimes t)v_{w\Lambda}=0.$$
If not, then $L(\Lambda)_{w\Lambda-\alpha_i+\delta}\ne 0$ and hence we find that $L(\Lambda)_{\Lambda-w^{-1}\alpha_i+\delta}\ne 0$. Now, $$(w\Lambda,\alpha_i)=m=(m\Lambda_0,w^{-1}\alpha_i)$$ implies that $w^{-1}\alpha_i=\alpha+\delta$ for some $\alpha\in R$
and hence $L(\Lambda)_{\Lambda-\alpha}\ne 0$.
 This forces $$\alpha\in R^+,\ \
L(\Lambda)_{s_\alpha(\Lambda-\alpha)}\ne 0. $$ Since $\Lambda=m\Lambda_0$ ,we find that $s_\alpha(\Lambda-\alpha)=\Lambda+\alpha$ and hence  $L(\Lambda)_{\Lambda+\alpha}\ne 0$
which is a contradiction.
This proves that $D(w\Lambda)$ is a quotient of $KR(m\omega_i)$. The fact that it as an isomorphism of $\lie g[t]$--modules is immediate from \cite[Theorem 1]{FoL}.

\end{pf}

\subsection{The case of exceptional Lie algebras} For the exceptional Lie algebras the definition and elementary properties of the modules $KR(m\omega_i)$ and their twisted analogs are exactly the same as for the classical Lie algebras. Moreover, the $\lie g$--module decompositions of the modules $KR(m\omega_i)$ (resp. $\krsm$) can be shown to coincide with the conjectural decompositions given in \cite{HKOTY}, \cite{HKOTT}, \cite{Kl}
as long as $i$ is such that the maximum value of $\varepsilon_i(\alpha)\le \check d_i$ for all $\alpha\in R^+$.

In the other cases, the main  difficulty lies in proving Proposition \ref{unupb}
(resp. Proposition \ref{twupb}). One reason for this is that the multiplicity of the irreducible representations occurring in a given graded component can be bigger than one.
The construction of the fundamental Kirillov--Reshetikhin modules
is also much more complicated since the number of irreducible components is large, in the case of $E_8$ and the trivalent node for instance, the total  number of components is conjectured to be three hundred and sixty eight with twenty four non--isomorphic isotypical components. In fact in this case, a general conjecture for the structure of the non--fundamental modules is not available and the combinatorics seems formidable.

\bibliographystyle{amsplain}

\end{document}